\documentstyle[12pt]{article}

\catcode `\@=11
\@addtoreset{equation}{section}

\catcode `\@=12

\newcommand{\be}{\begin{equation}}
\newcommand{\en}{\end{equation}}
\newcommand{\bea}{\begin{eqnarray}}
\newcommand{\ena}{\end{eqnarray}}
\newcommand{\beano}{\begin{eqnarray*}}
\newcommand{\enano}{\end{eqnarray*}}
\newcommand{\bee}{\begin{enumerate}}
\newcommand{\ene}{\end{enumerate}}

\newcommand{\R}{R \!\!\!\! R}
\newcommand{\N}{N \!\!\!\!\! N}
\newcommand{\Z}{Z \!\!\!\!\! Z}
\newcommand{\Hil}{{\cal H}}

\newcommand{\Vj}{{\cal V}_j}
\newcommand{\V}{{\cal V}}
\newcommand{\Wj}{{\cal W}_j}
\newcommand{\W}{{\cal W}}
\newcommand{\D}{{\cal D}}
\newcommand{\E}{{\cal E}}

\begin{document}
 
\thispagestyle{empty}

\vspace*{2cm}

    \begin{center}
  {\Large \bf Multi-Resolution Analysis in Arbitrary  \vspace{2mm}\\
   Hilbert Spaces} 
  \vspace{1.5cm}\\
   {\large F. Bagarello}
\vspace{1cm}\\
   Dipartimento di Matematica ed Applicazioni, Facolt\'a di Ingegneria,\\
   Universit\'a di Palermo, Viale delle Scienze, I-90128  Palermo, Italy.\\
E-mail:  bagarello@ipamat.math.unipa.it

   \end{center}

       \vspace{0.5cm}

  \begin{abstract}
  \noindent
We discuss the possibility of introducing a multi-resolution in a 
Hilbert space which is not necessarily a space of functions. We investigate which
of the classical properties can be translated to this more general framework and
the way in which this can be done. We comment on the procedure proposed by means of
many examples.   
\end{abstract} 

\vspace{3cm}


\vfill

\newpage

\section{Introduction}

The wavelet transform (WT) is by now a well established tool in many
branches
of physics, such as acoustics, spectroscopy, geophysics, astrophysics,
fluid
mechanics (turbulence),  medical imagery, \ldots (see 
\cite{wav3} for a survey of the present status). Basically it is a
time-scale
representation, which allows a fine analysis of nonstationary signals and a
good reconstruction of a signal from its WT, both in one and in two
dimensions
(image processing).

The basic formula for the (continuous) WT of a one-dimensional signal $s
\in
L^2(\R)$ reads: 
\be
 S(a,b) = 
     a^{-1/2} \; \int \; \overline{\psi}\left(\frac{x-b}{a}\right) \;
s(x)\;dx,
\label{cwt}
\en                             
where $a>0$ is a scale parameter and $b \in \R$ a translation  parameter.
Both the function $\psi(x)$, called the {\em analyzing wavelet},  and its
Fourier
transform $\hat\psi(\omega)$ must be well localized, and, in addition, 
$\psi$ is
assumed to have zero mean:
\be
\int \;\psi(x) \;dx = 0.
\en
 Combined with the localization properties, this relation makes the WT
 (\ref{cwt}) into a local filter and ensures its
efficiency  in signal analysis and reconstruction.

However, in practice, one often uses a {\em discretized} WT, obtained by
restricting the parameters $a$ and $b$ in (\ref{cwt}) to the points of a
lattice, typically a dyadic one:
\be
 S_{j,k} =  2^{-j/2} \; \int \; \overline{\psi}(2^{-j}x-k) \; s(x)\;dx,  
                          \quad j,k \in \Z.
\en
Very general functions $\psi(x)$ satisfying the admissibility conditions
described
above will yield a good WT, but then the functions
$\{ \psi_{j,k}(x) \equiv 2^{j/2}\psi(2^{j}x-k), j,k \in \Z\}$ are, in
general,
not orthogonal to each other! One of the successes of the WT was the
discovery
that it is possible to construct functions  $\psi(x)$ for which 
$\{ \psi_{jk}(x), j,k \in \Z\}$ is indeed an orthonormal basis of $ L^2(\R)$.
In addition, such a basis  still has the good properties of wavelets,
including space {\em and} frequency localization. This is the key to their
usefulness in many applications. In the past years a general way in 
which this function $\psi$ can be built up has been proposed, see \cite{dau,chui}
and \cite{mal} and references therein. This procedure is now known as a {\em
multiresolution analysis} (MRA) of $L^2(\R)$.

A MRA of
$L^2(\R)$ is an increasing sequence of closed subspaces $\{V_{j}\}_{j\in\Z}$
of $L^2(\R)$ such that
\begin{itemize}
\item[(1)] 
 $  \hspace{4mm} \ldots \subset V_{2} \subset V_{1} \subset V_0 \subset V_{-1}
\subset V_{-2}
  \subset \ldots;$
\item[(2)] 
  $\hspace{4mm}\bigcup_{j \in \Z} V_j \mbox{ is dense in } L^2(\R),$ 
 and $ \bigcap_{j \in \Z} V_j = \{0\};$
\item[(3)] 
$\hspace{4mm}f(x) \in  V_j \Leftrightarrow f(2x) \in  V_{j-1};$
\item[(4)] 
$\hspace{4mm}f(x) \in  V_0 \Rightarrow f(x-n) \in  V_{0}, \: \forall n\in \Z$;
\item[(5)]
\hspace{4mm}there exists a function $\phi \in V_0$, called a {\em scaling}
function,  such
that  $\{\phi(x-k), k \in {\Z}\}$ is an o.n. basis of $V_0$.
\end{itemize}

Each $V_j$ can be interpreted as an approximation space:  the approximation
of $f \in L^2(\R)$ at the resolution $2^{-j}$ is defined by its
projection onto $V_j$. The additional details needed for increasing the
resolution from  $2^{-j}$ to $2^{-(j+1)}$ are given by the projection of $f$
onto the orthogonal complement $W_j$ of $V_j$ in $ V_{j-1}$:
\be
V_j \oplus W_j =  V_{j-1},
\en
and we have:
\be
\bigoplus_{j \in {\Z}} W_j = L^2(\R).
\en
Then the theory asserts the existence of a function $\psi$, called the 
{\em mother} of the wavelets, explicitly computable from $\phi$, such that 
$\{\psi_{j,k}(x) \equiv 2^{j/2} \psi(2^jx-k), j,k \in \Z\}$ 
constitutes an orthonormal basis of $ L^2(\R)$: these are the 
{\em orthonormal wavelets}. 

The construction of  $\psi$ proceeds as follows. First, the inclusion $V_0
\subset V_{-1}$ yields the relation
\be   
\phi(x) = \sqrt{2} \,\sum_{n=-\infty}^\infty \; h_n \phi(2x - n),
\quad h_n = \langle \phi_{1,n} , \phi \rangle,
\label{11}
\en
which is known as "two-scale relation" (TSR). Taking its Fourier transforms, this
gives \be
\widehat\phi(\omega) = m_o(\omega/2) \widehat\phi(\omega/2),   
\label{12}
\en
where
\be
m_o(\omega) = \frac{1}{\sqrt{2}} \sum_{n=-\infty}^\infty \; h_n e^{-i n
\omega}
\en
is a $ 2\pi$--periodic function, belonging to $L^2([0,2\pi])$. Iterating
(\ref{12}), one gets the scaling function as the (convergent!) infinite product
\be
\widehat\phi(\omega) = (2\pi)^{-1/2} \prod_{j=1}^\infty \;
m_o(2^{-j}\omega).
\en
Then one defines the function $\psi \in W_0 \subset V_{-1}$ by the relation
\be
\widehat\psi(\omega) =  e^{i \omega/2} \; \overline{m_o(\omega/2 + \pi)}
                                          \;  \widehat\phi(\omega/2),
\en
or, equivalently,
\be   \label{psi}
\psi(x) = \sqrt{2} \; \sum_{n=-\infty}^\infty \;(-1)^{n-1} h_{-n-1} 
              \phi(2x - n),
\en
and proves that the function $\psi(x)$ indeed generates via dilations and
translations an o.n. basis with all the required properties.

Far from being only a mathematical tool, MRA has also many applications to physics; in
particular, for instance, it has been used by J.-P. Antoine together with the
author to approach the basis problem for the Fractional Quantum Hall Effect,
\cite{bag1}. Other applications of MRA to quantum mechanics are contained, for
instance, in \cite{jpa}, which also contains other interesting references.

As one may imagine already from the construction sketched above, all the proofs
of the Propositions concerning MRA are strongly related to the nature of the
objects we are dealing with. In particular, since we are considering
functions belonging to  $L^2(\R)$, we have at our disposal a very powerful tool:
the Fourier transform (FT). Indeed, if we take a look to reference \cite{dau}, we
see that FT is used to prove almost each step of the path leaving from
the MRA of $L^2(\R)$ to the required mother wavelet $\psi(x)$. Of course, the idea
of generalizing this structure to Hilbert spaces different from $L^2(\R)$, which,
moreover, are not necessarily function spaces,  forces the loss of this instrument.
The main content of this paper is the proof that, nevertheless, this generalization
can be performed in a quite natural and simple way, very much related to what we will
often call here the "classical" MRA, that is the MRA of $L^2(\R)$, and that many
results still hold true in the new setting. 

Of course, since the tools at our disposal are not so powerful as the ones in
$L^2(\R)$, we expect that some extra conditions will appear in the
statements of our Propositions.
We will discuss these conditions and we will show many examples in which they are
satisfied. 

\vspace{3mm}

The paper is organized as follows:

in the next Section we introduce the definition of multi-resolution (MR) of a
Hilbert space $\Hil$ and we give some examples.

In Section 3 we discuss which conditions on a given vector $\Phi \in \Hil$
implies that $\Phi$ generates a MR of $\Hil$, and the way in which this MR is
obtained.

In Section 4 we show how a MR produces a 'mother wavelet vector' and, therefore,
an o.n. basis of $\Hil$. We focus our attention to the new hypotheses which
need to be imposed in order to obtain this basis, considering some examples.

Finally, Section 5 is devoted to the discussion of different remarks concerning, for
instance, the construction of a scaling vector in $\Hil$, the approximation
scheme naturally carried by a MR and so on.

\section{Definitions and Examples}

In this Section we will introduce the definition of MR of a generic Hilbert
space $\Hil$, and we will show some examples fitting into this structure.

Let $\Hil$ be a (complex) Hilbert space and $\Vj$, $j\in \Z$, closed subspaces
of $\Hil$.

\noindent
	{\bf Definition 1}.--
We say that the set $\{\Vj\}_{j\in \Z}$ defines a MR of $\Hil$ if the following
properties are satisfied: 
\begin{itemize}
\item[(p1)]
 $ \hspace{4mm} \ldots \subset \V_{2} \subset \V_{1} \subset \V_0 \subset
\V_{-1} \subset \V_{-2}
  \subset \ldots.$
\item[(p2)] 
  $\hspace{4mm} \bigcup_{j \in \Z} \Vj \mbox{ is dense in } \Hil, $
and $ \bigcap_{j \in \Z} \Vj = \{0\}.$
\item[(p3)] 
$\hspace{4mm} \exists$ an unitary operator $\pi$ such that $$\varphi \in  \Vj
\Leftrightarrow \pi^j\varphi \in  \V_0 \quad \forall j \in \Z.$$ 
\item[(p4)] 
$\hspace{4mm} \exists$ an unitary operator $\tau$ such that $$\varphi \in  \V_0
\Rightarrow \tau^k\varphi \in  \V_0, \quad \forall k \in \Z.$$
Moreover we require that $\tau \pi = \pi \tau^2$. 
\item[(p5)]
\hspace{4mm} There exists a vector $\Phi \in \V_0$, called a {\em scaling} vector, 
such
that  $\{\tau^k \Phi, k \in {\Z}\}$ is an o.n. basis of $\V_0$.
\end{itemize}

\vspace{3mm}

{\bf Remarks}:-

(1) Of course the operators $\pi$ and $\tau$ play respectively the role of the
dilation and the translation in $L^2(\R)$. This is the reason why they have been
required to be unitary. For this same reason we have imposed the commutation law
in (p4). Sometimes, it will be useful to consider an extension of
this law, which can be proven by induction straightforwardly: for any $k\geq 0$
and for all $l\in \Z$ we have: \be
\left\{
       \begin{array}{ll}
        \tau^l\pi^k = \pi^k \tau^{2^kl}, \\
        \tau^{2^kl}\pi^{-k} = \pi^{-k} \tau^{l}.
       \end{array}
       \right.
\label{21}
\en

(2) The definition above can be modified, as in the canonical case, leaving
unchanged its meaning. In particular we can substitute property (p5) above with
the following equivalent condition:

\vspace{3mm}

(p5$'$) \hspace{4mm} there exists a vector $\varphi \in \V_0$  such
that  $\{\tau^k \varphi, k \in {\Z}\}$ is a Riesz basis of $\V_0$.

\vspace{3mm}

This condition, thought being apparently weaker than (p5), is completely
equivalent to this one. In fact, any o.n. basis is obviously a Riesz basis.
Conversely, using the procedure extensively discussed in \cite{mey} and true in
any Hilbert space, starting from the vector $\varphi$ in (p5$'$) we can construct a
vector $\Phi \in \V_0$ satisfying condition (p5).

\vspace{3mm}

Of course a mathematical definition makes sense if there exists some non-trivial
structure which fits all the requirements. Therefore, before going on with
considering the consequences of Definition 1, we discuss some examples. 

\vspace{3mm}

{\bf Example 1 }-

The first obvious example is obtained going back to the canonical situation. We
take $\Hil = L^2(\R)$, $\tau$  the translation operator (indicated with
$T$) defined as $(\tau f)(x)=(Tf)(x)=f(x-1)$, $\pi$ the dilation operator
(indicated with $P$) defined by $(\pi f)(x)=(Pf)(x)=\sqrt{2}f(2x)$ and $\Vj=V_j$,
$\{V_j\}$ being a MRA of $L^2(\R)$. Finally, the scaling vector $\Phi$ coincides
with the scaling function $\Phi(x)$ of the canonical definition.

\vspace{3mm}

{\bf Example 2 }-

The second example is obtained considering an unitary
map $U$  between the Hilbert spaces $L^2(\R)$ and $\Hil$. We can
define  
\be
\Vj=U V_j, \quad \tau=UTU^{-1},\quad \pi=UPU^{-1},
\label{21bis}
\en
where $T$, $P$ and $V_j$ have just been introduced in Example 1, and we take the
scaling vector $\Phi$ to be the $U$-image of the scaling function
$\Phi(x)$, $\Phi= U\Phi(x)$.

With these positions it is quite easy to verify that all the requirements of
the Definition 1 are satisfied. 

Incidentally we observe that this situation is really quite general and, in a
sense, it solves the question of the existence of a MR in spaces different from
$L^2(\R)$: in fact, since we are always considering separable Hilbert spaces,
$\Hil$ is necessarily isomorphic to $L^2(\R)$ so that an
unitary map between the two spaces can be defined without difficulties: it only
needs to map an o.n. basis of $L^2(\R)$ into an o.n. basis of $\Hil$. 

We may wonder now if Definition 1 implies something interesting in the
analysis of (separable) Hilbert spaces. We believe that this is so,
the reason being that, how we will show in Section 4, the existence of
a MR in a Hilbert space $\Hil$ implies the existence of an o.n.
basis of $\Hil$ which can be found explicitly starting from the MR  itself, without
going trough any unitary map between $\Hil$ and $L^2(\R)$.  In other words, if
$\Psi_{j,k}(x)$ is an o.n. basis of $L^2(\R)$ constructed starting with a mother
wavelet $\Psi(x)$ and $U$ is an unitary map between  $L^2(\R)$ and $\Hil$, then
$U\Psi_{j,k}(x)$ is an o.n. basis of $\Hil$ which can be written as
$\pi^j\tau^k(U\Psi(x))$, $\tau$ and $\pi$ being defined in (\ref{21bis}). In Section
4 we will learn how to built up the same o.n. basis of $\Hil$, that is essentially
the vector $U\Psi(x)\in \Hil$, using as the only ingredient the MR of $\Hil$, without
considering at all the action of the operator $U$. This can be
convenient, at least in some situations: if, for instance, we consider the Fock
space ${\cal F}$ of an infinite number of bosons, it is not evident at all which
should be the unitary map between ${\cal F}$ and $L^2(\R)$. On the contrary, it
may happen that we are able to build up a MR in this space using, for instance,
the procedure discussed in the next Section, and, therefore, to get directly an o.n.
basis of ${\cal F}$.

\vspace{3mm}

{\bf Pre-example 3 }-

Let us suppose that between the Hilbert spaces $L^2(\R)$ and $\Hil$ it exists
an invertible map $U$ which is bounded together with its inverse. With the same
definitions of the previous example, see (\ref{21bis}), we obtain that the spaces
$\{\Vj\}$ satisfy almost all the properties of Definition 1. We observe, by the
way, that:

- the o.n. basis $\{T^k \Phi(x)\}$ of $V_0$ is mapped in a Riesz basis of
$\V_O$;

- the operators $\tau$ and $\pi$ are not unitary.

While the first point does not imply any serious consequence, the second
does. The reason, which will appear clear in the next Section, is that we cannot
restrict our search for an o.n. basis of the whole $\Hil$ only to $\W_0$, that
is to the orthogonal complement of $\V_0$ into $\V_{-1}$. This is why we consider
this one only as a {\em pre-example}. The natural question is therefore the following:
given the spaces $\Vj=UV_j$, can we define two unitary maps $\tilde \pi$ and
$\tilde \tau$ different from $\pi$ and $\tau$, acting on $\Vj$ as in Definition 1?
We will briefly comment on this point in the last Section, where we will prove
the existence of a pair of unitary operators in a given Hilbert space, which
satisfies the commutation property $\tau \pi=\pi \tau^2$ and which are not
unitary equivalent to $T$ and $P$.

\vspace{3mm}

{\bf Example 4 }-

Let us consider an o.n. basis $\{\varphi_n(x)\}$, $n\in \N_0$, in $L^2(\R)$ (e.g.
the Hermite polynomials). In correspondence we can define the following set of
coefficients: \be
b^{(j,k)}_l\equiv 2^{-j/2}\int^{2^j(k+1)}_{2^jk}\overline{\varphi_l(x)}\; dx.
\label{22}
\en
It is easy to show, first of all, that any sequence $b^{(j,k)}$ belongs to
$l^2(\N_0)$. More explicitly we have:
$$
<b^{(j,k)},b^{(j,k')}>=\delta_{k,k'}.
$$
We introduce now two operators $\pi$ and $\tau$, acting on these sequences,
defined by
\be
\pi b^{(j,k)} = b^{(j-1,k)}, \hspace{1cm} \tau b^{(j,k)} = b^{(j,k+2^{-j})},
\quad \quad \forall j,k \in \Z,
\label{23}
\en
and we put $\Phi \equiv b^{(0,0)}$. The vector $b^{(j,k)}$ satisfies the
following TSR: 
\be
b^{(j,k)}= \frac{1}{\sqrt{2}}\left(b^{(j-1,2k)}+b^{(j-1,2k+1)}\right),
\quad \quad \forall j,k \in \Z.
\label{24}
\en
It
is also straightforward to verify that, if we define for any $j\in \Z$ the
following subsets of $l^2(\N_0)$, $$
\Vj \equiv \mbox{linear span} \overline{\{b^{(j,k)}:\:\: k\in \Z\}},
$$
then the set $\{\Vj\}$ defines a MR of $l^2(\N_0)$. 

Incidentally we observe that this example may be considered as a particular
case of Example 2, where $U$ is the unitary map between the o.n. bases 
$\{\varphi_n(x)\}$ and $\{e_n\}$ (the canonical basis in $l^2(\N_0)$),
and is nothing but the {\em image} in $l^2(\N_0)$ of the Haar MR of 
$L^2(\R)$. 

We also notice that, with the same techniques, we can construct examples in
spaces different from $l^2(\N_0)$ like, for instance, in the Fock 
space of one boson, $\Hil_{boson}$. This can be done quite simply by considering
the natural isomorphism $I$ between $l^2(\N_0)$ and $\Hil_{boson}$, defined
in the following canonical way: let $\Psi_0$ be the ground state of the theory and
$a^\dagger$ the creation operator. Then, any $\Psi\in \Hil_{boson}$ can be
expanded in terms of the o.n. basis $\Psi_k\equiv
\frac{(a^\dagger)^k}{\sqrt{k!}}\Psi_0$ as $\Psi=\sum_k c_k \Psi_k$. We conclude
 defining $I(\Psi)=\{c_k\}_{\{k\in \N_0\}}$, which, of course, belongs to
$l^2(\N_0)$.

\vspace{3mm}

{\bf Example 5 }-

The last example we want to discuss is, perhaps, the most interesting. The
Hilbert space in which we want to construct a MR is $L^2(I)$, where $I=]0,1[$.
Therefore, we are considering an Hilbert space of functions, but we are no
longer working on the whole real line.         

Let $v(x)$ be a differentiable function defined in $I$ with value in the whole
$\R$, strictly monotone. Obviously the inverse function $v^{-1}(x)$ exists,
and it maps $\R$ into $I$. For concreteness sake we choose
$$
v(x) = \log \left(\frac{x}{1-x}\right), \quad \quad v^{-1}(x)=\frac{1}{1+e^{-x}}.
$$  
By means of the function $v^{-1}(x)$ we define subsets of $I$ in the following
way:
\be
S_{j,k} \equiv [v^{-1}(2^jk),v^{-1}(2^j(k+1))[.
\label{25}
\en
Furthermore, we consider the normalized function in $L^2(I)$:
\be
\rho(x) = \left\{
         \begin{array}{ll}
                 \frac{1}{\sqrt{x(1-x)}} \quad x\in S_{0,0}\\
                 0 \quad \mbox{ otherwise,}
       \end{array}
        \right. 
\label{26}
\en
and we define the following maps on $L^2(I)$:
\be
\tau \varphi(x) \equiv \varphi(v^{-1}(v(x)-1))\frac{\sqrt{e}}{x(1-e)+e} = 
\varphi\left(\frac{x}{x(1-e)+e}\right)\frac{\sqrt{e}}{x(1-e)+e},
\label{27}
\en
\be
\pi \varphi(x) \equiv \sqrt{2}\,
\varphi(v^{-1}(2v(x)))\frac{\sqrt{x(1-x)}}{x^2+(1-x)^2} = 
\sqrt{2}\,
\varphi\left(\frac{x^2}{x^2+(1-x)^2}\right)\frac{\sqrt{x(1-x)}}{x^2+(1-x)^2}, 
\label{28} \en
for all $\varphi \in L^2(I)$. It is an easy computation to verify that both $\tau$
and $\pi$ are unitary in $L^2(I)$ and that they satisfy the commutation rule $\tau
\pi=\pi \tau^2$.

The reason why we have introduced the function $\rho$ is that this behaves nicely
under the action of the operators $\tau$ and $\pi$:
\be
\rho_{-l,k}(x) := (\pi^{-l}\tau^k\rho)(x) = \left\{
         \begin{array}{ll}
                 \frac{2^{-l/2}}{\sqrt{x(1-x)}} \quad x\in S_{l,k}\\
                 0 \quad \mbox{ otherwise.}
       \end{array}
        \right. 
\label{29}
\en
Making use of the fact that $S_{j,k}\bigcap S_{j,k'}=\emptyset$ for all $k$
different from $k'$, equation (\ref{29}) implies that the set 
$\{(\tau^k\rho)(x),\: k\in \Z\}$ is made of o.n. functions. Moreover $\rho(x)$
satisfies the following TSR
\be
\rho(x) = \frac{1}{\sqrt{2}}\left(\pi \rho(x)+\pi\tau \rho(x)\right),
\label{210}
\en
which looks very much the same as the TSR satisfied by the characteristic
function of the interval $[0,1[$, appearing in the Haar MRA of $L^2(\R)$. At this
point, if we define the following sets: $$
\Vj \equiv \mbox{linear span} \overline{\{\pi^{-j}\tau^k\rho(x) :\:\: k\in \Z\}},
$$
we can prove that $\{\Vj\}$ defines a MR of $L^2(I)$. Almost all the requirements are
obvious; the only point which need to be discussed is property (p2). The
easiest way to show its validity is the following: since the support $S_{j,k}$ goes
to zero when $j\rightarrow -\infty$ then $\cap_{j\in \Z} \Vj = \{0\}$. Moreover,
the density of the set $\cup_{j\in \Z} \Vj$ in $L^2(I)$ follows again from the 
properties of $S_{j,k}$ which, when $j$ and $k$ change, 'cover' the whole $I$
with intervals whose measure varies with $j$ and $k$. As a matter of fact,
we are replacing the density in  $L^2(\R)$ of the functions piecewise constant
with the density in $L^2(\R)$ of the functions which are "piecewise 
$\frac{1}{\sqrt{x(1-x)}}$".

\section{Scaling Vector as Starting Point}

In reference \cite{dau} it is discussed the way in which a given function
$\Phi(x)$ belonging to $L^2(\R)$ which satisfies a TSR 
\be 
\Phi(x) = \sqrt{2} \sum_n h_n \Phi(2x-n),
\label{31}
\en
where $\sum_n|h_n|^2<\infty$, and for which there exist two positive constants
$0<\alpha\leq\beta<\infty$ such that  
$$
\alpha\leq \sum_{n\in \Z}|\tilde\Phi(\omega+2\pi n)|^2\leq \beta,
$$
can be used to build up a MRA of $L^2(\R)$. We begin with defining the spaces
$V_j$ as the closure of the linear span of the vectors
$P^{-j}T^k\Phi(x)$, for $k$ ranging in $\Z$. Here $P$ and $T$ are the dilation
and translation operators introduced in the Example 1 of the previous Section. It
is easy to show that almost all the properties of a MRA of $L^2(\R)$ are
automatically verified following this construction, whenever  $\Phi(x)$ satisfies the
two conditions above. The only property which needs some extra care is the density
of $\cup_j V_j$ in $L^2(\R)$.  Nevertheless, this is a consequence of an extra
condition on $\Phi(x)$, namely the condition of being $\tilde \Phi(\omega)$
bounded for all $\omega$ and continuous near $\omega=0$, with  $\tilde \Phi(0)\neq
0$.

\vspace{3mm}

For an arbitrary Hilbert space $\Hil$, an analogous procedure can still be
considered. We consider as a starting point a scaling vector $\Phi
\in \Hil$ which satisfies a certain TSR \be 
\Phi = \sum_{n\in \Z} h_n \pi\tau^n\Phi,
\label{32}
\en
with $\sum_{n\in \Z}|h_n|^2 <\infty$, and such that the vectors $\tau^n\Phi$,
$n\in \Z$, are mutually orthonormal. Here $\tau$ and $\pi$ are \underbar{any} given
pair of unitary operators such that \be
\tau \pi = \pi \tau^2.
\label{33}
\en
As we have already discussed, the orthonormality condition can be weakened requiring
that the set $\{\tau^n\Phi\}$ is a Riesz set in a given subspace of $\Hil$.

Following the same steps as in the Example 5 of Section 2 we define the
following closed subsets of $\Hil$: 
\be
\Vj \equiv \mbox{linear span} \overline{\{\pi^{-j}\tau^k\Phi :\:\: k\in \Z\}}.
\label{34}
\en
It is easy to see that these subsets satisfy conditions (p1), (p3), (p4) and (p5)
of the Definition 1. It remains to understand in which hypotheses also condition
(p2) is satisfied. We discuss now two Propositions which, together, imply (p2).

\vspace{3mm}

\noindent
	{\bf Proposition 1}.--
Let $\D$ be a dense subset of $\Hil$. Let us assume that corresponding to the
scaling vector $\Phi$ there exist:

\noindent
- a function $m(x)$ going to zero for $x\rightarrow \infty$;

\noindent
- a sequence $\{A_{f,k}\}$ summable and depending on the arbitrary vector $f\in
\D$, such that 
\be
|<\pi^jf,\tau^k\Phi>|^2\leq A_{f,k} \,m(j), \quad \forall j,k\in \Z.
\label{35}
\en
Therefore we have $\cap_{j\in \Z}\Vj =\{0\}$.

\vspace{5mm}

\noindent
	{\underline {Proof}}

Calling $P_j$ the projection operator on $\Vj$ we know that, for any $f\in \D$
$$
\|P_jf\|^2=\sum_{k\in\Z}|<\pi^jf,\tau^k\Phi>|^2 \leq m(j) \sum_{k\in\Z} A_{f,k}
\rightarrow 0,
$$
for $j$ going to infinity. This implies that, when $j
\rightarrow \infty$, $P_j\rightarrow 0$  strongly on $\Hil$,
which, together with condition (p1) of Definition 1, implies the statement.
\hfill $\Box$

\vspace{3mm}

{\bf Example 6 }-

We take $\Hil=L^2(\R)$, $\D={\cal S}(\R)$ and $\Phi(x)$ the usual characteristic
function in $[0,1[$. It is easily seen that $$
|<\pi^jf,\tau^k\Phi>| \leq 2^{-j/2} \int_{2^jk}^{2^j(k+1)}|f(x)|\,dx.
$$
At this point, with a minor generalization of the hypotheses of Proposition 1, we
put 
$$
A_{f,k}^{(j)} \equiv \left( \int_{2^jk}^{2^j(k+1)}|f(x)|\,dx \right)^2, \quad
m(x)\equiv 2^{-x}.
$$
We observe that the dependence on $j$ in $A_{f,k}^{(j)}$ really plays no
role, since it disappear after the sum over $k$ is performed. In fact we have, for
any function $f\in {\cal S}(\R)$ 
$$
\|P_jf\|^2\leq 2^{-j} \sum_{k\in\Z}\left( \int_{2^jk}^{2^j(k+1)}|f(x)|\,dx
\right)^2 \leq 2^{-j} \left( \sum_{k\in\Z} \int_{2^jk}^{2^j(k+1)}|f(x)|\,dx
\right)^2 =  2^{-j} \|f\|_1^2,
$$
which is finite since $f$ belongs to  ${\cal S}(\R)$. Of course, the right hand
side goes to zero when $j$ diverges to $+\infty$, so that $\cap_{j\in \Z}\Vj =
\{0\}$.

\vspace{3mm}

In reference \cite{micch} an interesting proof of this fact is contained for
$\Hil=L^2(\R)$, which seems to be, in a certain sense, less related to the nature of
the Hilbert space. By the way, even its generalization to generical Hilbert spaces
appear to be not so straightforward.

The next Proposition gives a sufficient condition for the union of the $\Vj$ to
be dense in $\Hil$.

\vspace{3mm}

\noindent
	{\bf Proposition 2}.--
Let us suppose that the scaling vector $\Phi$ satisfies the following property:

\noindent
-for any $f\in \Hil$
there exist two integers $j_0$ and $k_0$ and a strictly positive constant $c$ such
that \be
|<\pi^{j_0}f,\tau^{k_0}\Phi>|\geq c \|f\|.
\label{36}
\en
Then $\cup_{j\in \Z} \Vj$ is dense in $\Hil$.

\vspace{5mm}

\noindent
	{\underline {Proof}}

In order to prove the statement we will show that the orthogonal
complement of the set $(\cup_{j\in \Z} \Vj)$ in $\Hil$ contains only the zero
vector. 

Let $f\in \Hil\backslash (\cup_{j\in \Z} \Vj)$. This means that $f\notin \Vj$ for
any choice of $j$ in $\Z$. Therefore, since the set
$\{\pi^{-j}\tau^k\Phi\}_{\{k\in \Z\}}$ is an o.n. basis of $\Vj$,
$$
<f,\pi^{-j}\tau^{k}\Phi>=0, \quad \quad \forall j,k \in \Z.
$$
Using the hypothesis of the Proposition we conclude that
$$
\|f\| \leq \frac{1}{c}|<\pi^{j_0}f,\tau^{k_0}\Phi>|=0,
$$
which implies that $f=0$.
\hfill$\Box$

\vspace{3mm}

The above results suggest the following obvious consideration: what is (almost) for
free in $L^2(\R)$ seems to have a price in $\Hil$! What we can say at the moment is
that, even if
the hypotheses of Propositions 1 and 2 appear rather strong, expecially the one of
the second Proposition, we hope to approach succesfully the problem of weakening
these conditions in a future paper.

\section{Some Results}

As we have already anticipated, one of the main utility of the definition of a MR
of a Hilbert space $\Hil$ is that it allows to build up an o.n. basis in $\Hil$,
exactly as it does in $L^2(\R)$.
We will devote this Section to show how the existence of this basis can be
proved. This proof, of course, requires techniques which are completely different
from the classical ones, which strongly rely on the Fourier transform, \cite{dau}. We
will show that the same results as in the classical case can be found also in this
more abstract situation, at least if some extra conditions are required. We will also
show, however, that these new conditions are rather weak and satisfied in many
relevant examples.

The first step consists in defining new subspaces of $\Hil$: we call $\Wj$
the orthogonal complement of $\Vj$ in $\V_{j-1}$:
\be
\Wj \oplus \Vj =  \V_{j-1}.
\label{41}
\en
For these spaces we can state the

\vspace{3mm}

\noindent
	{\bf Proposition 3}.--
The spaces $\Wj$ satisfy the following properties:

\begin{itemize}
\item[(a)]
 $\Wj$ is orthogonal to $\W_{j'}$ for any $j$ different from $j'$;
\item[(b)]
the direct sum $\oplus_{j\in \Z} \Wj$ is dense in $\Hil$;
\item[(c)]
$\varphi \in \Wj \iff \pi^j\varphi \in \W_0, \quad \forall j\in \Z$;
\item[(d)]
$\varphi \in \W_0 \iff \tau^k\varphi \in \W_0, \quad \forall k\in \Z$.
\hfill$\Box$

\end{itemize}

\vspace{4mm}

The proof of this Proposition is quite similar to the canonical one, so that we
will omit it here. 

The consequences of this Proposition are the same as in the canonical case. In
particular, if a set $\{\tau^k \Psi\}_{\{k\in \Z\}}$ is an o.n. basis of $\W_0$
conditions (b) and (c), together with the unitarity of $\pi$, ensure us that the
set $\{\pi^j \tau^k  \Psi\}_{\{j,k\in \Z\}}$ is an o.n. basis of the whole Hilbert
space. It is therefore enough to find an o.n. basis in
$\W_0$, again as in the canonical situation, in order to obtain an o.n. basis of
$\Hil$. It is worthwhile to observe that it is really at this stage that the
unitarity of $\pi$ is required, and this is the reason why the pre-example 3
cannot be considered as a {\em real} example.

Let us now remind that, since by construction the set $\{\pi
\tau^n\Phi\}$ is an o.n. basis in $\V_{-1}$ and $\Phi \in \V_0 \subset \V_{-1}$,
we can write for $\Phi$ the usual two-scale relation
\be
\Phi = \sum_{n\in \Z} h_n \pi \tau^n\Phi,
\label{42}
\en
where the coefficients $h_n$ are 
\be
h_n = <\Phi, \pi \tau^n\Phi>.
\label{43}
\en
We define now the infinite matrix $K$ whose elements are
\be
    K_{l,n}  = \left\{
          \begin{array}{ll}
            h_{l-n},    &       \mbox{if } n \mbox{ is even in } \Z, \\
           (-1)^{l-1}h_{-l+n},   &       \mbox{if } n \mbox{ is odd in } \Z. 
       \end{array}
        \right. 
\label{44}
\en
The role of $K$ will be clear in the following Proposition.

We are now ready to discuss one of the main results of this paper, which shows
how any MR of a given Hilbert space is related to an o.n. basis of $\Hil$, at
least under some conditions which are usually satisfied.

We have the following

\vspace{3mm}

\noindent
	{\bf Proposition 4}.--
Let  $\{\Vj\}_{j\in \Z}$ be a MR of $\Hil$ as in Definition 1. In particular
let $\Phi$ be the scaling vector satisfying the two-scale relation (\ref{42}).
Defining the (mother wavelet) vector
\be
\Psi \equiv \sum_{n\in \Z}(-1)^{n-1} h_{-n-1}\pi \tau^n\Phi,
\label{45}
\en
then:
\begin{itemize}
\item[(a)]
 $\{\tau^k \Psi\}_{k\in \Z}$ is a set of mutually orthonormal vectors;
\item[(b)]
 if $h_n \in \R$ then $\Psi \in \W_0$;
\item[(c)]
 if $h_n \in \R$ and if the matrix $K$ is invertible then 
$\{\tau^k \Psi\}_{\{k\in \Z\}}$ is an o.n. basis of $\W_0$.
\end{itemize}

\vspace{5mm}

\noindent
	{\underline {Proof}}

\begin{itemize}
\item[(a)]
The first statement follows from the orthonormality requirement on the scaling
vector $\Phi$. It is easy to obtain, see \cite{dau,chui} for $\Hil=L^2(\R)$, that
\be
<\Phi, \tau^k \Phi> = \sum_{n\in \Z}h_n\overline{h_{n-2k}}=\delta_{k,0}.
\label{46}
\en
On the other hand, we get
\be
<\Psi, \tau^k \Psi> = \sum_{n\in \Z}h_n\overline{h_{n+2k}},
\label{47}
\en 
which implies, together with (\ref{46}), the orthonormality of the set 
$\{\tau^k \Psi\}_{\{k\in \Z\}}$.

\item[(b)]
In order to prove that $\Psi$ belongs to $\W_0$ we have to show that $\Psi$
belongs to $\V_{-1}$ and that it is orthogonal to all the vectors of the form
$\tau^k \Phi$. The first requirement is obvious, since $\Psi$ is defined exactly in
terms of an o.n. basis of $\V_{-1}$. For what concerns the second point we 
notice that
$$
<\Psi, \tau^k \Phi>= \sum_{n\in \Z}(-1)^{n-1}h_{-n-2l-1}\overline{h_{n}}=2i
\Im\left( \sum_{n\in \Z}h_{2n}\overline{h_{-2n-2l-1}} \right),
$$
where $\Im(z)$ indicates the immaginary part of the complex number $z$. The statement
follows immediately.

\item[(c)]
We only have to prove that the set $\{\tau^k\Psi, \:\: k\in\Z\}$ is complete in
$\W_0$. This can be proved showing that any vector $\chi \in \W_0$ orthogonal to
$\tau^k\Psi$ for all integer $k$ is necessarily zero.

Since $\W_0 \subset \V_{-1}$ we can expand $\chi$ in terms of the o.n. basis of
$\V_{-1}$, $\{\pi \tau^n\Phi\}$: 
$$
\chi = \sum_{n\in \Z} a_n \pi \tau^n\Phi, \quad \quad a_n=<\chi, \pi \tau^n\Phi>.
$$
Moreover, since $\chi$ belongs to $\W_0$, it is necessarily orthogonal to all the
vectors $\tau^n\Phi$, for any $n\in \Z$. And yet, due to our hypothesis, $\chi$
is also orthogonal to all $\tau^n\Psi$, for any $n\in \Z$. Using the reality of
the coefficients $h_n$, these conditions imply that:
$$
<\chi,\tau^n\Psi> = 0 \quad \Rightarrow \quad \sum_l (-1)^{l-1}a_lh_{-l+m}=0,
$$
for all odd integers $m$, and 
$$
<\chi,\tau^n\Phi> = 0 \quad \Rightarrow \quad \sum_l a_lh_{l-m}=0,
$$
for all even integer $m$. Recalling the definition (\ref{44}) of the matrix $K$ we
obtain the condition
$$
\sum_l a_lK_{l,m}=0, \quad m\in \Z,
$$
which, of course, implies that all the coefficients $a_l$ are zero if, and only if,
the matrix $K^{-1}$ does exist. 
\hfill$\Box$
\end{itemize}

\vspace{3mm}

{\bf Remarks}:-

(1) As we have already discussed, the vector $\Psi$ in (\ref{45}) is such that the
set $\{\pi^j\tau^k\Psi:\quad j,k\in \Z\}$ is an o.n. basis in $\Hil$.

(2) It is interesting to observe that the reality condition on the coefficients
$h_n$ really need to be imposed, even if in the classical case it was not
necessary. In fact we have the following simple counterexample:

let our scaling vector be such that all the coefficients $h_n$ are zero in the
TSR but for $n=0,1,2,3$. In particular, we take $\overline{h}_0=h_2
=\frac{1}{\sqrt{4}}i$ and $h_1=h_3=\frac{1}{\sqrt{4}}$. With this choice the
orthonormality condition for the vectors $\tau^k
\Psi$, $\sum_nh_n\overline{h}_{n-2k} = \delta_{k,0}$ is satisfied.  Nevertheless
not all the coefficients are real and, as a consequence, $\Psi$ is not orthogonal to
all the vectors $\tau^k\Phi$! In particular, for instance, it is easy to see that
$<\Psi, \tau^{-1}\Phi>=-\frac{i}{2}\neq 0$.

Of course the reason for this 'pathological' behavior relies in the lack of
the extra conditions on the coefficients $h_n$, $\sum_n h_{2n}= \sum_n
h_{2n+1}=1$, see \cite{dau,chui}, which in our context there is no reason to
introduce, and which are not satisfied by the example above. 

By the way we observe that in all the examples in $L^2(\R)$, at least in our
knowledge, the TSR contains only real coefficients.

(3) Other solutions are possible for the mother wavelet vector above, as in the
classical case. One possibility, which will be  used in the examples below, is the
following: \be
\tilde \Psi \equiv \sum_{n\in \Z}(-1)^{n} h_{-n+1}\pi \tau^n\Phi.
\label{48}
\en
For $\tilde \Psi$ we can prove all the same results of Proposition 4 more
or less  by means of the same steps.

\vspace{3mm}

{\bf Example 4 (reprise)}-

The above theorem allows us to write down an o.n. basis in $l^2(\N_0)$. The
procedure is all contained in equation (\ref{48}) and in the TSR
(\ref{24}) satisfied by the vectors $b^{(j,k)}$: both of these equations implies
that the mother wavelet vector is
$$
\tilde \Psi = \frac{1}{\sqrt{2}}\left(b^{(-1,0)}-b^{(-1,1)}\right),
$$
so that the o.n. basis in $l^2(\N_0)$ is given by the following set $\{\frac{1}
{\sqrt{2}} \left(b^{(-1-j,2k)}-b^{(-1-j,2k+1)}\right),\:\: j,k \in \Z\}$.

\vspace{3mm}

{\bf Example 5 (reprise)}-

Recalling the TSR in (\ref{210}) and equation (\ref{45}) we can write down
the mother wavelet function even for this problem:
\be
\Psi_\rho(x)  = \left\{
         \begin{array}{ll}
                 \frac{1}{\sqrt{x(1-x)}} \quad x\in S_{-1,-1},\\
                 -\frac{1}{\sqrt{x(1-x)}} \quad x\in S_{-1,-2},\\
                 0 \quad \mbox{ otherwise.}
       \end{array}
        \right. 
\label{49}
\en
Of course another mother wavelet can be obtained by means of equation (\ref{48})
and is deduced by $\Psi_\rho(x)$ simply replacing $S_{-1,-1}$ with $S_{-1,0}$
and  $S_{-1,-2}$ with $S_{-1,1}$.

\section{Final Remarks and Conclusions}

We begin this last Section with the following remark: also in a general Hilbert
space a MR can be useful to obtain controlled approximations of any given vector of
$\Hil$. In fact, by means of the projections operators $P_j$ on the spaces $\Vj$, and
of the other  projections operators $E_j\equiv P_{j-1}-P_j$, which project on $\W_j$,
we see that for any $\varphi \in \Hil$ we can write the following equality:
$$
\varphi-P_j\varphi = \oplus_{k=j}^{-\infty}E_k\varphi.
$$
Introducing the coefficients $c_{k,l}\equiv <\varphi,\pi^{-k}\tau^l\Psi>$, where
$\Psi$ is the mother wavelet vector in (\ref{45}), we see that
$$
\|\varphi-P_j\varphi\|^2=\sum_{k=j}^{-\infty}\sum_{l\in \Z}|c_{k,l}|^2,
$$
which, of course, can be made as small as we want simply by taking $-j$ big
enough. This implies that any vector in $\Hil$ can be approximated with the
precision required simply considering its projection on an opportune space
$\Vj$. This also implies that, if the vector we want to approximate already
belongs to a certain $\Vj$, its approximation coincides, from a certain point on,
with the vector itself.

\vspace{3mm}

The second remark is the following: as for the classical MRA of $L^2(\R)$, also in
this abstract situation we can modify a bit some of the hypotheses of the
Definition 1. For instance, we can start with two unitary operators $\tilde \pi$
and $\tilde \tau$ which, instead of the condition in (p4), satisfy the following
commutation rule: $\tilde \tau \tilde \pi= \tilde \pi \tilde \tau^3$. Under this
condition in $L^2(\R)$ can be constructed a MR, the main difference being that, for
each $j\in \Z$, two orthogonal spaces $W_j^1$ and  $\W_j^2$ must be introduced such
that $V_{-1}=V_0\oplus W_0^1 \oplus  W_0^2$. We claim that the same procedure can be
used also for arbitrary Hilbert spaces. We believe that this situation can also be
generalized to all the pairs of unitary operators $(\hat\pi,\hat\tau)$,
with $\hat \tau \hat \pi= \hat \pi \hat \tau^k$, $k$ being a natural number bigger
than one.

\vspace{3mm}

As already announced previously in the paper, we are going to comment on the
existence of the pair of unitary operators $\pi$ and $\tau$ satisfying the
condition $\tau\pi=\pi\tau^2$. In particular we want to show that two such
operators need not to be the unitary transformations of $P$ and $T$, as in Example
2. To show this fact, it is enough to consider the following easy counterexample: let
$\Hil = C\!\!\!\!C^2$. Such an Hilbert space is, obviously, not unitarily
equivalent to $L^2(\R)$. Let us consider the following unitary operators on $\Hil$:
$$
\tau = e^{i\frac{2\pi}{3}\hat b\cdot \underline{\sigma}}, \quad \pi = i e^{ia_0}
\hat a \cdot \underline{\sigma}.
$$
Here $\hat a$ and $\hat b$ are two mutually orthogonal vectors, both normalized.
It is straightforward to see that, indeed, $\tau  \pi=  \pi  \tau^2$, though
they cannot be written as $\tau=UTU^{-1}$ and $\pi=UPU^{-1}$ for any unitary
operator $U:L^2(\R)\rightarrow C\!\!\!\!C^2$, simply becouse such an $U$ does not
exist.

Of course, the existence of a similar pair of operators is still open when
dim$(\Hil)=\infty$, which is the only relevant situation in which a MR can be
introduced. We plain to consider again this problem in a future paper.

\vspace{3mm}

The last point we want to discuss in this Section is the way in which, under
certain (heavy, indeed!) hypotheses on the operators $\tau$ and $\pi$, with $\tau 
\pi=  \pi  \tau^2$, we can construct a scaling vector $\Phi$ satisfying the TSR
(\ref{32}) and such that the vectors $\tau^n\Phi$, $n\in \Z$, are mutually
orthonormal.

This problem can be approached in a particularly simple way whenever the operator
$\tau$ admits an o.n. basis of eigenvectors $\varphi_l$:
\be
\tau \varphi_l = t_l \varphi_l.
\label{51}
\en
In this condition, assuming, as usual, that all the eigenvalues $t_l$ are
different from zero, we have
\be
\tau^k \varphi_l = t_l^k \varphi_l, \quad \forall k\in \Z.
\label{52}
\en
To get the solution, it is enough to find a vector $\Phi$ such that
\be
     \left\{
      \begin{array}{ll}
         <\Phi, \tau^k\Phi>=\delta_{k,0} \quad \forall k\geq 0, \\
          \Phi = \sum_{n\in \Z}h_n\pi\tau^n\Phi.
       \end{array}
        \right. 
\label{53}
\en
We expand $\Phi$ in terms of the eigenvectors of $\tau$:
\be
\Phi=\sum_lc_l\varphi_l, \hspace{3cm} c_l=<\Phi,\varphi_l>,
\label{54}
\en
and we substitute this expansion in (\ref{53}). In particular, the first equation
becomes
\be
\sum_l|c_l|^2t_l^k=\delta_{k,0}, \quad \forall k\geq 0,
\label{55}
\en
while the second can be written in the form
\be
c_k=\sum_lc_lq_l\pi_{l,k},
\label{56}
\en
where we have defined $\pi_{l,k}\equiv <\pi\varphi_l,\varphi_k>$ and $q_l\equiv
\sum_nh_nt_l^n$. (To avoid convergence problems we can suppose, at this stage,
that only a finite number of $h_n$ are different from zero.) Introducing
furthermore an infinite matrix $\lambda$ with matrix elements $\lambda_{k,l}:=
q_l\pi_{l,k}$, the matrix
\beano
\E=  \left(
\begin{array}{ccccccc}
1 & 1 & 1 & 1 & . & . & .  \\ 
t_1 & t_2 & t_3 & t_4 & . & . & .  \\ 
t_1^2 & t_2^2 & t_3^2 & t_4^2 & . & . & .   \\ 
t_1^3 & t_2^3 & t_3^3 & t_4^3 & . & . & .   \\
. & . & . & . & . & . & .   \\
. & . & . & . & . & . & .   \\
. & . & . & . & . & . & .   \\
\end{array}
\right)
\enano
and the vector $\gamma$ whose elements are the square moduli of the $c_l$,
$\gamma_l=|c_l|^2$, the above conditions can be rewritten in the following
matricial form: \be
\E \gamma = u
\label{57}
\en
where $u$ is the vector with only the first component equal to 1 and all the
others equal to zero, and
\be
c=\lambda c.
\label{58}
\en
We assume, at this stage, that the determinant of the matrix $\E$ can be computed
without convergence problems. Therefore, if all the eigenvalues $t_k$ are different,
 it is known that this determinant is different from zero, so that
$\E^{-1}$ exists. Therefore, the vector $\gamma$ can be obtained: $\gamma=\E^{-1}
u$. The scaling vector $\Phi$ can be found if, among all the solutions $c$ related to
this unique $\gamma$, it exists at least one vector which satisfies equation
(\ref{58}).

At our actual level of knowledge we have still no a priori argument ensuring
 that such a scaling vector $\Phi$ exists (if we don't want to make use of the
isomorphism between $\Hil$ and $L^2(\R)$ which certainly solves the existence
question!). What we have proposed here is only a way in which, if everything works,
such a vector can be obtained. By the way, the technique which we have used here to
approach the problem is not totally satisfactory, also because the classical situation
in $L^2(\R)$ does not verify the main hypothesis. The reason is simply that the
translation operator does not admit an o.n. family of functions belonging to
$L^2(\R)$. For all these reasons we believe that this construction procedure for
$\Phi$ really requires an extra effort, to be completely satisfactory. We hope to be
able to consider this problem again in more details in a future paper.

\vspace{6mm}

We have seen how a MR can be introduced for general Hilbert spaces, not
necessarily of functions, like, for instance, $l^2(\N_0)$ and $L^2(0,1)$. Our
procedure appears to be an economical generalization of the classical MRA of spaces
$L^2(\R)$. Obviously, some differences arise somewhere, since in general we have
in mind to deal with abstract vectors which may not be functions. 
Among the other topics, we have discussed the possibility of obtaining an o.n.
basis of $\Hil$ related to a given MR. We have also discussed the way
in which a MR can be built up starting with an opportunely chosen (scaling) vector of
$\Hil$. As we have already said, a lot of things can still be done: Propositions 1
and 2 both contain very strong hypotheses which, in our opinion, could be weakened.
In the same way, Proposition 4 requires conditions on $h_n$ and $K$ which should
be better understood. Moreover, we plain to undertake in a future paper a deeper
analysis of the construction of the scaling vector $\Phi$, which we have only
sketched here, and of all the pairs of unitary operators $\pi$ and $\tau$ such that
$\tau\pi=\pi\tau^2$ in an infinite dimensional Hilbert space. We also plain to
translate into this framework many of the
 features already well established in the classical situation.

\vspace{50pt}

\noindent{\large \bf Acknowledgments} \vspace{5mm}

The author wishes to thank Dr. C. Trapani and Prof. H. Fujita for their useful
suggestions. He is also indebted with M.U.R.S.T. for financial support.

\end{document}